# An alternative to modules. *Une alternative aux modules.* [*][†]

Gabriele Ricci

*via dei Fiori 11, CH–6600 Muralto, Switzerland.*

e–mail: gabrielericci@bluewin.ch

**Abstract**

Recently, a geometrical characterization of vector spaces served to generalize them into a new class of algebras. Instead of the *algebraic* properties of the underlying fields, we generalized the recently discovered property *of such spaces* that generates these fields. This property also concerned the semi-linear transformations, which are necessary to define geometrical invariance.

Yet, the class of such geometrical generalizations was practically unknown. We only knew that it differs from the class of modules, because of a simple example (the sum monoid of natural numbers).

Here, we partly clarify the extent of this class: we prove that it contains at least the one of all commutative based universal algebras and we provide it with four new examples.

Again, our further examples are not modules. They exhibit a wider choice of both types and algebraic properties, though they keep the representations by "coordinates".

Typeset by LaTeX

## 0 Preliminaries

**0.0 Introduction**. If we define vector spaces as two-components systems, we usually generalize them as modules by generalizing their latter component, a field, into a ring. Yet, a recent characterization of vector spaces, which we recall in 3.4 (A), allows us a one-component definition as mere universal algebras with a simple geometric property. (Thus, the old huge list of equations that links the two components and defines the field disappears [16, 17].)

In fact, the new definition exploits the notion of (general) dilatations for universal algebras, which generalize the dilatations we know from vector spaces, as we recall in 2.0. This reminds us of the "generalized conception of space" that A.N. Whitehead proposed in the introduction of [22] as a tool for investigating all universal algebras (including the ones that generalize such spaces).

General dilatations replace the old postulates for the field with its construction: just as the field was seen to underlie a vector space, this one-component view now sees *the vector space as underlying its field*. Besides, they are *fundamental in all Universal Algebra:* without them (e.g. by mere automorphisms) one cannot even get invariance, as [18, 19] proved by generalizing the well-known geometric counterexample in [21] from a vector space on the complex field, i.e. from a two-components algebra, to a new (one-component) universal algebra.

This construction is not peculiar to vector spaces. We can use it on certain other algebras to get a new generalization of fields: the "endowed dilatation monoids" that 2.2 will recall. Therefore, the algebras that generate such "endowed monoids", called dilatation full algebras, do generalize vector spaces, yet by Whitehead's geometric approach and in a

---

closer agreement with the Erlanger Program of F. Klein. Then, we can consider them a case of *geometrical generalizations* of vector spaces.

Nevertheless, the extent of dilatation full algebras is unknown. From the example in 1.8 of [17] we only know that it is not the one of modules.

To narrow this gap, we will prove that *at least all commutative based universal algebras are dilatation full.* (Our proof needs some set-theoretical extensions of usual algebraic notions.) Four such algebras will show how they generate wide-ranging "endowed monoids".

**0.1 Conventional notation.** We consider functional composition as a restriction of relational composition, here denoted by $\cdot$, namely $f \cdot g$ is the composition "of $g$ and $f$" and

$$(0) \qquad (f \cdot g)(x) = f(g(x)).$$

This does not agrees with the definition **4.7** of [8], our reference for Set Theory, nor with *the ones of some textbooks* of Algebra that see $f \cdot g$ as the composition "of $f$ and $g$".

We write $f:A \to B$ to say that $f$ is a function with arguments in the whole set $A$ and values in $B$, $f:A \rightarrowtail B$ or $f:A \twoheadrightarrow B$ to say that it also is one to one or onto $B$ and $f:A \rightarrowtail\!\!\!\twoheadrightarrow B$ to say it is a bijection onto $B$. Often, we will replace "function" with "indexing", to emphasize values (while avoiding "family"). $\mathbb{I}_S$ denotes the identity function on $S$. Then, $\mathbb{I}$ denotes the function (a proper class) that provides every set $S$ with its identity function.

For sets $A$ and $X$, $A^X$ will denote the set of functions from $X$ to $A$, although when $X$ is a natural number $n = \{0, \ldots, n-1\}$ we might well take $A^n$ to be the set of $n$-tuples on $A$. A similar ambiguity will concern $A^{n' \times n''}$, a set of two-indices arrays.

**0.2 Combinatory notation.** Contrary to the previous notation $A^n$, the two following notations will not confuse set-theoretically different objects. They will locally omit some set-theoretical specification that still can be formally recovered from their context. (Boldfacing symbols will remind us of this.)

(A) Among the functions in $A^X$ we consider the constant ones. In case $A \neq \emptyset$ we denote the constant function with a value $a \in A$ by $\boldsymbol{k}_a$:

$$(1) \qquad \boldsymbol{k}_a(x) = a, \text{ for all } x \in X \neq \emptyset.$$

Yet, (1) defines a function $\boldsymbol{k}:A \to A^X$ even in all other cases: for $X = \emptyset$ and $A \neq \emptyset$ there are the trivial cases $\boldsymbol{k}_a = \emptyset$, while for $A = \emptyset$ we get $\boldsymbol{k} = \emptyset$. We call $\boldsymbol{k}$ the *constant generating function (for $A$ and $X$).*

We will use constant generating functions for several different $A$ and $X$. As typefaces run out, we will not distinguish their notation. E.g., given sets $Y$, $B$ and $a \in A$, we write the identities

$$(2) \qquad \boldsymbol{k}_a \cdot M = \boldsymbol{k}_a, \qquad \text{for all } M:X \to Y \text{ and}$$
$$(3) \qquad M \cdot \boldsymbol{k}_a = \boldsymbol{k}_{M(a)}, \quad \text{for all } M:A \to B,$$

which follow either immediately, when $X = \emptyset$, or from (1), when $X \neq \emptyset$. Yet, our quantifications distinguish three functions with the same notation: in (2) $\boldsymbol{k}:A \to A^Y$ on the left, whereas $\boldsymbol{k}:A \to A^X$ on the right, while in (3) $\boldsymbol{k}:A \to A^X$ on the left, whereas $\boldsymbol{k}:B \to B^X$ on the right.

(B) Another incomplete yet useful notation serves to generalize the medial law $f(g(m_{0,0}, m_{0,1}, \ldots m_{0,s-1}), g(m_{1,0}, m_{1,1}, \ldots m_{1,s-1}), \ldots g(m_{r-1,0}, m_{r-1,1}, \ldots m_{r-1,s-1})) = g(f(m_{0,0},$ $m_{1,0}, \ldots m_{r-1,0}), f(m_{0,1}, m_{1,1}, \ldots m_{r-1,1}), \ldots f(m_{0,s-1}, m_{1,s-1}, \ldots m_{r-1,s-1}))$ for all natural numbers $r$ and $s$ and all $m{:}r \times s \to A$, which defines commutativity for (finitary) operations $f{:}A^r \to A$ and $g{:}A^s \to A$ as in **1.1.5** of [9], to the case of $f{:}A^R \to A$ and $g{:}A^S \to A$, where the ranks $R$ and $S$ have to be arbitrary sets.

In fact, for the new (generalized) $f$ and $g$ we can introduce a new $m{:}R \times S \to A$, but we can no longer write this law as before. Yet, now we can replace this $m$ with an $m'{:}R \to A^S$, such that $m'_i(j) = m_{i,j}$ for all $j \in S$ and $i \in R$, to rewrite the left hand side of the medial law as $f(g \cdot m')$ by (0). Similarly for the right hand side: we start from an $\tilde{m}{:}S \times R \to A$, with $\tilde{m}_{j,i} = m_{i,j}$ to define another $m''{:}S \to A^R$ by $m''_j(i) = \tilde{m}_{j,i}$ for all such $j$ and $i$. Then, the required generalization is the identity $f(g \cdot m') = g(f \cdot m'')$ for all $m$. Note that, when $R$ and $S$ are natural numbers $r$ and $s$, we get the previous medial law by the conventions in `0.1`.

This generalization exploits the transition from $m'{:}R \to A^S$ to $m''{:}S \to A^R$ that permutes the applications in $m'_i(j) = m''_j(i)$, which is similar to an array transposition. In fact, it corresponds to an operator that swaps rows and columns in $m{:}R \times S \to A$ to get $\tilde{m}{:}S \times R \to A$. However, as we are dealing with several functions ($m$, $\tilde{m}$, $m'$ and $m''$), a specific notation might be convenient.

Further, in `0.4 (A)` this transposition-like transition will serve to generalize Menger's superposition [2], which concerns finitary operations. In fact, here the finiteness restriction on arities will fail even for finitary algebras: from their finitary operations we will derive other operations that might need infinite arguments. (See also **2.5** in [20].) It will avoid finiteness restrictions also in `0.6 (E)`. Thus we make the following definition.

Given any $m{:}I \to A^J$, $\mathbf{C}_m^{(J)}$ denotes the *exchanged function of* $m$, $\mathbf{C}_m^{(J)}{:}J \to A^I$, defined by $[\mathbf{C}_m^{(J)}(j)]_i = m_i(j)$ for all $j \in J$ and $i \in I$. When $I \neq \emptyset$, $m$ determines $J$. Then, we can simplify this notation as $\boldsymbol{c}_m : J \to A^I$, which, when $I$ and $J$ are fixed, defines an *exchange function* $\boldsymbol{c}{:}(A^J)^I \to (A^I)^J$. However, we will use several different exchange functions, again without distinguishing their notation.

Moreover, we will also leave the duty of specifying $J$ to the context, when $I = \emptyset$. Then, all such conventions allow us to write

(4) $\qquad \mathbf{C}_\emptyset^{(J)} = \boldsymbol{c}_\emptyset = \begin{cases} \emptyset & \text{when } J = \emptyset \text{ and} \\ \boldsymbol{k}_\emptyset{:}J \to 1 & \text{when } J \neq \emptyset, \end{cases}$

and in general

(5) $\qquad [\boldsymbol{c}_m(j)]_i = m_i(j) \quad \text{for all } j \in J \text{ and } i \in I.$

In addition to typeface savings, such conventions highlight functional features of algebraic interest better than a set-theoretically complete notation. For instance, they exhibit how two (different) $\boldsymbol{c}$ work the same (by (5)) in the identity

(6) $\qquad \boldsymbol{c}_{\boldsymbol{c}(m)} = m, \quad \text{for all such } m,$

which follows from (5) and implies $\boldsymbol{c} \cdot \boldsymbol{c} = \mathbb{I}_{(A^J)^I}$, and

(7) $\qquad \boldsymbol{c}{:}(A^J)^I \mapsto\!\!\!\twoheadrightarrow (A^I)^J.$

In spite of this functional notation, our $\boldsymbol{c}$ and $\boldsymbol{k}$ are not (set-theoretical) functions, but define them by the context. When we implicitly redefine such symbols, we exploit the exchange combinator ("elementary permutator" in [1]) and the constant generator ("elementary canceller" ibid.) within a set-theoretical setting. (Often, combinators show a working sameness even stronger than the one we observed in (6), in the sense that it concerns functions with more different algebraic or categorical properties. This makes combinators a heuristic aid to find these functions and properties in spite of such differences.)

**0.3 Algebras**. On a carrier set $A$ we consider operations, which we call *set-ary,* where a rank $S$ of an operation $f{:}A^S \to A$ can be any set. When we have to consider conventional finitary operations, we assume that they are replaced with them. E.g., an $f{:}A^2 \to A$ will replace an $f'{:}A \times A \to A$ by the natural map for $A \times A \simeq A^2$, while $f$ might keep the possible infix notation of $f'$. We will deal with both the ("fundamental") operations of an algebra and others. When we do not declare that an operation is of the algebra, we consider both cases.

The algebras of our main concern here do not really need an indexing of their operations. Such an indexing serves for a category of algebras with *homomorphisms,* whereas now we mainly consider single (universal) algebras with endomorphisms. Still, as this indexing occurs in an auxiliary algebra (see 0.6 (E)), we conform to the indexed case.

(Our key definitions and the main theorem in 2.6 will exploit this auxiliary algebra. Yet, one cannot formalize it by conventional operations, even when only conventional operations define the algebras of our main concern, see 0.7 (C). Moreover, even some simpler algebras of our concern might need unconventional operations, e.g. the complete union semilattices on all subsets of infinite sets.)

Therefore, here we define an *algebra on $A$* as an indexing $\alpha \in \prod_{i \in I} A^{A^{r(i)}}$ of such operations, where $r$, its algebra type, is an arbitrary function. By $\mathcal{O}$ we denote its set of operations: $\alpha{:}I \twoheadrightarrow \mathcal{O}$. Yet, contrary to conventional definitions, $\mathcal{O}$ will determine $A$. Thus, if all operations in $\mathcal{O}$ are nullary, $A$ cannot be larger than the set of their values, which might happen by defining our algebra as a pair $\langle A, \alpha \rangle$.

We assume $\mathcal{O} \neq \emptyset$, since an $\mathcal{O} = \emptyset$ gives the class of all sets as its carrier, according to some definitions of a carrier (when such a definition has the prefix " for all $f \in \mathcal{O}$"). Since the case $A = \emptyset$ mainly concerns the initial settings in computer implementations, we will allow uninterested readers to skip it by putting the observations relevant to it into square brackets. As this is the only case where an operation (or algebra) does not determine its rank (or type), outside such brackets "*a* rank/type of" will become "*the* rank/type of".

Our result as well as some of the intermediate statements concern *based algebras,* which are defined also by one of their possible bases in addition to an indexing of operations. As 0.8 will show, our bases, which we will define in 0.6, are equivalent to the free generating "families" of § 24 in [5] for the case of finitary algebras.

**0.4 Definitions.** (A) We define *the (set-ary) composition of an indexing $G{:}S \to A^{A^Y}$ of $Y$-ranked operations with an operation $g{:}A^S \to A$* as the function $\ell = g \cdot \boldsymbol{c}_G$. Then, by (7) $\ell{:}A^Y \to A$ and $\ell(M) = g(\boldsymbol{c}_G(M))$ for all $M{:}Y \to A$, where, for a nullary $g$, by (3) and (4)

(8) $$S = \emptyset \quad \text{implies} \quad \ell = \boldsymbol{k}_{g(\emptyset)}.$$

[Notice that, when $A = \emptyset$, $g$ cannot be nullary, while $G$ requires $Y \neq \emptyset$.] Clearly, when $S$ and $Y$ are natural numbers $s$ and $y$, the natural bijections for $A^S \simeq A^s$ and $A^Y \simeq A^y$ allow us to rewrite our composition as a (finitary) superposition. Notice that we are not requiring either $g \in \mathcal{O}$ nor $G{:}S \to \mathcal{O}$.

(B) Given an algebra on $A$ and a set $Y$, $\mathcal{L}_Y$ will denote the set of the *Y-ary elementary functions* of the algebra, defined as the functions we get by such compositions with its operations from all the *projections* $p_x:A^Y \to A$ defined by

(9) $\qquad\qquad p_x(M) = M(x) \quad \text{for all } x \in Y \text{ and all } M:Y \to A \ .$

Formally, $\mathcal{L}_Y = \bigcap \{\mathcal{F} \subseteq A^{A^Y} \mid p:Y \to \mathcal{F} \text{ and } (f \in \mathcal{O}, f:A^S \to A \text{ and } L:S \to \mathcal{F} \text{ imply } f \cdot \boldsymbol{c}_L \in \mathcal{F})\} \subseteq A^{A^Y}$. When the arity $Y$ is a natural number, $\mathcal{L}_Y$ corresponds to a *proper* subset of the algebra clone ([2], **VI** of [3]).

When $Y = \emptyset$, the only compositions involved are the ones with an indexing of nullary constants from the empty set of projections. Without any nullary $g$ this gives $\mathcal{L}_\emptyset = \emptyset$. In general, $\mathcal{L}_\emptyset$ is the set of nullary constants that corresponds to the subalgebra closure of the empty set. [When $A = \emptyset$ and $Y \neq \emptyset$, all projections are empty and $\mathcal{L}_Y = \{\emptyset\}$.]

(C) If we equalize all the arguments in any $Y$-ary elementary function for $Y \neq \emptyset$, we get a function in $A^A$. Formally, the composition of $\boldsymbol{k}:A \to A^Y$ and any $\ell \in \mathcal{L}_Y$ gives us a function $\ell' = \ell \cdot \boldsymbol{k}:A \to A$ and this defines a function $j:\mathcal{L}_Y \to A^A$ such that $j_\ell = \ell'$.

Here, we can forget $Y$, since from (B) we easily see that we get all such $\ell'$ for all $Y \neq \emptyset$ by merely setting $Y = 1$. Then, $\mathcal{L}' \subseteq A^A$ will denote the set of all such $\ell'$. Moreover, $Y = 1$ makes $\boldsymbol{k}$ into a bijection $\boldsymbol{k}:A \rightarrowtail\!\!\!\twoheadrightarrow A^1$ that makes $j$ too into a bijection $j:\mathcal{L}_1 \rightarrowtail\!\!\!\twoheadrightarrow \mathcal{L}'$. We call any $\ell' \in \mathcal{L}'$ a *rank-less elementary function*.

(D) As in the finitary case, we say that an indexing $U:X \to A$ is an *(indexed) generator* of our algebra $\alpha$, when every $a \in A$ is the value $\ell(U)$ for some $X$-ary elementary function $\ell \in \mathcal{L}_X$. Also, the independence of any $U:X \to A$ (viz. $\ell(U)$ determines the whole $\ell$ for all $\ell \in \mathcal{L}_X$) is to require *an elementary function generator,* defined as a function $\chi:A \twoheadrightarrow \mathcal{L}_X$ such that $\chi_{\ell(U)}(M) = \ell(M)$ for all such $\ell$ and $M:X \to A$. Clearly, an independent $U$ is a generator (which determines a conventional basis) if and only if there exists a single $\chi$.

0.5 Analytic representations. A part of Linear Algebra, which one century ago was called "Analytic Geometry", concerns the calculus of the "usual" vector space matrices (two-indices arrays). If among them we consider the square ones, this calculus provides the endomorphism monoid of a vector space and other related structures with convenient representations.

Unfortunately, when the matrix notion originated, neither the Theory of Data Structures nor Universal Algebra were born. Then, it was possible to think of a matrix (a mathematical object enjoying certain properties in vector spaces) as a case of an array (a simple data structure). However, our square matrices enjoy all their properties (including the one of possibly being two-indices arrays as in 0.7 (A)) only because they represent all endomorphisms in a certain way.

Still, many present textbooks of Linear Algebra introduce such matrices directly as (square) arrays, while endomorphism representations, if mentioned, come later. Anyway, such a matrix notion is not general: in many based universal algebras one is not able to represent endomorphisms as such arrays. (See 0.4 of [14] for several examples outside vector spaces that differ from two-indices arrays but still are able to model well-known objects.)

Endomorphism representations are not peculiar to vector spaces. Here, every our algebra (hence every universal algebra) will have them, provided only that it has a basis as such a vector space does. (In a sense these representations also concern the universal algebras that lack bases, because the ones that can have them are all and only the free ones.) Yet, since we cannot call them "linear", we call them "analytic" representations, while we keep the word "matrix" to denote a representation result.

However, one should not expect to keep the two-indices: in general the single index that will serve the basis will be enough. (We only will keep the optional adjective "square" as a possible reminder that we do not consider homomorphisms, which within vector spaces conventionally have rectangular arrays.) Since bases can have any cardinality, such an index might be neither finite nor denumerable.

Conversely, as 0.8 will show, if a generator provides an algebra with such an endomorphism representation *by the same functional construction* used in a vector space, then it is an (indexed) basis according to the conventional definitions recalled in 0.3 and 0.4 (D). This construction will generalize the starting idea "matrix = two-indices array" even in vector spaces, where it comes from. In fact, as detailed in 0.7 (B), in some vector spaces our generalized matrices, which become "one-index arrays", represent endomorphisms in a more natural way.

Therefore, when we use analytic representations for universal algebras, we can conveniently define bases by this endomorphism representability.

**0.6 Definitions**. Let $\mathbb{E}_\alpha \subseteq A^A$ be the set of all endomorphisms of an algebra $\alpha$ on a set $A$. Given a set $X$, let $U{:}X \to A$ and consider the function $\boldsymbol{r}_U{:}\mathbb{E}_\alpha \to A^X$, defined by $\boldsymbol{r}_U(h) = h \cdot U$, for $h \in \mathbb{E}_\alpha$. Namely $\boldsymbol{r}_U$ "samples each $h$ at" $U$ by providing each $x \in X$ with the value $h(U(x))$. A generator $U{:}X \to A$ used to define such a sampling will be called a *frame of $\alpha$*.

If this sampling serves to represent every endomorphism by any sample and conversely, namely if our algebra satisfies
$$\boldsymbol{r}_U{:}\mathbb{E}_\alpha \rightarrowtail\!\!\!\twoheadrightarrow A^X, \tag{10}$$
then every structure on $\mathbb{E}_\alpha$ defines another on $A^X$. While the former structure is abstract, the latter will depend on frame $U$.

(As shown in 0.7 (B), within vector spaces this dependence conflicts with the treatments that uniformly define matrices only as two-indices arrays even when the vectors are not one-index arrays. Yet, it will not conflict with our more complete treatment.)

In particular, we consider endomorphism application: for each $a \in A$, its application to any $h \in \mathbb{E}_\alpha$, is given by the function $q{:}A \to A^{\mathbb{E}_\alpha}$ such that $q_a(h) = h(a)$ for all such $a$ and $h$. Then, $\boldsymbol{r}_U$ has to represent $q$ by a function $\chi{:}A \to A^{A^X}$ such that $\chi_a(M) = q_a(h)$ whenever $\boldsymbol{r}_U(h) = M$. The construction of such a $\chi$ will again involve the exchange combinator, which in 0.8 will give us back the single $\chi$ of 0.4 (D).

Therefore, *if a frame $U$ satisfies (10),* then we say that

**(A)** $\boldsymbol{r}_U$ is an *analytic representation of $\mathbb{E}_\alpha$*, while $X$ is its *dimension set* and the cardinality of $X$ is its *dimension* (see below),

**(B)** its inverse $\eta = \boldsymbol{r}_U^{-1}{:}A^X \rightarrowtail\!\!\!\twoheadrightarrow \mathbb{E}_\alpha$, which extends any sample assignment $M \cdot U^{-1}$ onto the endomorphism $h = \eta_M$ with $h(U_x) = M_x$ for all $x \in X$, is the *(sample) extension function from $U$*,

**(C)** $A^X$ is the set of the *(general square) matrices of $\alpha$ with respect to $U$*, while every value $M(x)$ of a matrix $M{:}X \to A$ is its *column at $x \in X$*,

**(D)** $U$ is a *basis* or *(general) reference frame* of $\alpha$, while its columns $U(x)$ are *reference elements*, which form the *basis set* $B \subseteq A$ for $U{:}X \twoheadrightarrow B$,

**(E)** the *algebra of the conjugate functions derived from* $\alpha$ with respect to $U$ is the indexing $\chi{:}A \to A^{A^X}$, defined by (10) from the previous function $q$ as $\chi_a(\boldsymbol{r}_U(h)) = h(a)$ for $h \in \mathbb{E}_\alpha \subseteq A^A$ and $a \in A$, namely by **(B)**

(11) $\qquad \chi_a(M) = \eta_M(a), \quad$ for all $M{:}X \to A$ and $a \in A$,

which by (5) implies $\chi = \boldsymbol{c}_\eta$,

**(F)** while the value $\chi_a{:}A^X \to A$ of this indexing at an $a \in A$, which is an operation as in **0.3** of this constant-type algebra, is the *function conjugate of $a$ with respect to $U$.*

[Notice that $A = \emptyset$ by (10) implies $X = \emptyset$,] whereas for a singleton $A$ every set $X$ satisfies (10). In the former case we say that *the carrier (of the algebra) is trivial;* in the latter that *the algebra is trivial.* By (10), when the algebra is not trivial, $X = \emptyset$ if and only if $\mathbb{E}_\alpha = \{\mathbb{I}_A\}$. This happens when all algebra elements are constants. By (11) this also implies that $\chi{:}A \rightarrowtail A^1$ merely is the generator of singleton constants:

(12) $\qquad X = \emptyset \;$ iff $\; \chi_a(M) = a = \boldsymbol{k}_a(M), \;$ for all $M{:}X \to A, a \in A$.

In **(B)** notice that $M \cdot U^{-1}$ always is a function. In fact, when the algebra is trivial this (relational) composition gives a (singleton) function, while without triviality we still get a function, since $U^{-1}$ is a function, i.e. $U{:}X \rightarrowtail A$, because $U(x') = U(x'')$ implies that $[\boldsymbol{r}_U(h)]_{x'} = [\boldsymbol{r}_U(h)]_{x''}$ for all $h \in \mathbb{E}_\alpha$, i.e. $M_{x'} = M_{x''}$ for all $M \in A^X$, which implies $x' = x''$.

In **(A)** recall that, contrary to the case of a vector space which has a single dimension for all its representations, a non-trivial universal algebra allows several cases: no representation dimensions, when no basis exists; single representation dimension, either finite or not, and infinitely many finite dimensions [4, 6].

In the last case, one cannot always say that the algebra has any of such numbers as one of *its algebra dimensions.* In order to be relevant to the algebra, not just to the representation, such a number must not change under all transformations of the based algebra onto itself, which *can be more general* than the automorphisms. Such an invariance is not guaranteed: there are based algebras of the last kind both with and without (algebra) dimensions, see **3.6** respectively in [19] and [18]. (Then, these counterexamples *deny any "geometrical" significance to automorphism groups in Universal Algebra,* just as [21] did in Linear Algebra.)

**0.7 Examples.** **(A)** Take $A$ as the set of the usual $n$–tuples of elements of a field, $\alpha$ as their vector space "over" the same field and $X$ as $n$. In this case any frame $U{:}X \to A$ amounts to the selection of $n$ linearly independent vectors from $A$, which are $n$–tuples. Then, $U$, or any other $M{:}X \to A$, corresponds to an $n \times n$ usual matrix (an array) with columns $U_x$ or $M_x$.

If we use the frame $U$ that corresponds to the Kronecker matrix (viz. $U_0 = (1, 0, \ldots, 0)$, $U_1 = (0, 1, \ldots, 0)$, ..., $U_{n-1} = (0, 0, \ldots, 1)$), then for any endomorphism $h$ of the vector space and any $x \in X = n$ the endomorphic image $h(U_x)$ gives the $x$-th column vector $M_x = h(U_x)$ of the bijective representation of $h$ by a usual matrix. Hence, $U$ is a reference frame and by (11) $\chi_a(M)$ is the product of vector $a$ times the matrix $M$.

Therefore, the conjugate function $\chi_a$ of a vector $a$ is similar to its linear form. The only difference is that the former acts on vectors, while the latter acts on field numbers. (It follows that the conjugate functions in a vector space, as well as in any based universal algebra, in addition to their algebra in **0.6 (E)**, form another algebra that always is isomorphic to the

8                                                                                                    *G. Ricci*

starting one, as in 6.6 of [11], whereas linear forms merely form the adjoint space as in II.3 of [0].)

(B) On the contrary, when the representation of endomorphisms in (10) still concerns a based *vector space,* but with an arbitrary carrier and/or an arbitrary reference frame, it gives general matrices, possibly different from two-indices arrays, and different conjugate functions. Since this is less familiar than the former case, we will disregard it as an example for vector spaces and we will call *usual* the former vector space, as well as its corresponding structures. Yet, *even in vector spaces* to be a (general) matrix is not be a two-indices array, as reference vectors need not to be one-index arrays.

We can always replace two-indices arrays for such arbitrary matrices, because we can *transform* the latter (arbitrary) vector spaces into the former ones, not because the latter matrices do not exist *nor lack any use.* No more than fifty years ago, the inverse transformations (target vectors were functions implemented by electric waves with a proper spectrum) provided computing, also for Linear Algebra, with tools that were competitive with the digital computations of the time.

(We are referring to the "Analog Computing". To solve a numerical problem, engineers working in several computation laboratories preferred designing an electric circuit to programming digital computers. Putting plugs or even hard-wiring and hand-welding were faster and more reliable than exploiting the then available compilers. The idea of a matrix for such engineers might have been deeper than the one in the Linear Algebra textbooks that prefer to introduce a matrix as an array than as a homomorphism representation.)

(C) Another interesting case, which concerns vector spaces, is the infinite dimensional one. When $X$ is not finite, contrary to the case in (A), the conjugate function $\chi_a$ of a vector $a$ is not anymore similar to its linear form. Indeed, the latter always has to index its arguments by some *finite* subset $X'_a \subset X$, not by our $X$. Hence, while the vector space operations and its linear forms are conventional, the algebra of the conjugate functions derived from it fails to be conventional.

Notice also that replacing the finite ordinal $n$ of (A) with an infinite one (and not with a mere arbitrary infinite set) is an attempt to keep the order restriction, which might be misleading or useless even in the finite case. In fact, in a real life case the indexing of reference vectors could need a structure different from an order or no structure at all. While a finite ordinal at least serves to recall the order by which pens write arguments, in general such an "ink-theoretical" motivation might disappear.

0.8 Theorem. *Given a generator $U{:}X \to A$, there exists its (single) elementary function generator $\chi{:}A{\twoheadrightarrow}\mathcal{L}_X$ as in 0.4 (D) if and only if (10) holds. When $\chi$ exists, it is the algebra of the conjugate functions in 0.6 (E), which has the same set of endomorphisms as the original algebra: $\mathbb{E}_\chi = \mathbb{E}_\alpha$. Then,*

(13) $\qquad\qquad h \in \mathbb{E}_\alpha \quad \text{iff} \quad h(\chi_a(M)) = \chi_a(h \cdot M) \quad \text{for all } a \in A \text{ and } M{:}X \to A.$

*Proof.* (If) Take $\chi$ as in (11) and $h = \eta_M$ for each $M{:}X \to A$. Since $h \in \mathbb{E}_\alpha$, it also is an endomorphism of every $\ell \in \mathcal{L}_X$, see also 1.1. Hence, $h(\ell(M')) = \ell(h \cdot M')$ also for $M' = U$. Then, by (5) $\chi_{\ell(U)}(M) = \eta_M(\ell(U)) = h(\ell(U)) = \ell(h \cdot U) = \ell(M)$, which uniquely defines $\chi{:}A{\twoheadrightarrow}\mathcal{L}_X$ since $U$ is a generator.

(Only if) $\boldsymbol{r}_U$ is one to one. In fact, let $\eta = \boldsymbol{c}_\chi{:}A^X \to A^A$. Then, for all $a \in A$ and $h \in \mathbb{E}_\alpha$, by (5) $[(\eta \cdot \boldsymbol{r}_U)(h)]_a = \eta_{h \cdot U}(a) = \chi_a(h \cdot U) = \chi_{\ell(U)}(h \cdot U) = \ell(h \cdot U) = h(\ell(U)) = h(a) = [\mathbb{I}_{\mathbb{E}_\alpha}(h)]_a$, where $\ell \in \mathcal{L}_X$ and $\ell(U) = a$. It also is onto $A^X$: for all $x \in X$ and $M{:}X \to A$, $[(\boldsymbol{r}_U \cdot \eta)(M)]_x = [\boldsymbol{r}_U(\eta_M)]_x = \eta_M(U_x) = \chi_{U(x)}(M) = \chi_{p_x(U)}(M) = p_x(M) = M_x = [\mathbb{I}_{A^X}(M)]_x$ by (6), (9) and since $p_x \in \mathcal{L}_X$.

$\mathbb{E}_\alpha \subseteq \mathbb{E}_\chi$, since every $\chi_a$ is an $X$-ary elementary function. $\mathbb{E}_\alpha \supseteq \mathbb{E}_\chi$ by (10), since for all $a \in A$, if $h \in \mathbb{E}_\chi$, $h(a) = \eta_M(a)$ with $M = h \cdot U \in A^X$. In fact, $h(\chi_a(U)) = \chi_a(h \cdot U)$, viz. by (11) $h(\eta_U(a)) = \eta_{h \cdot U}(a)$, where $h(\eta_U(a)) = h(\eta_{\boldsymbol{r}_U(\mathbb{I}_A)}(a)) = h(\mathbb{I}_A(a)) = h(a)$.   Q.E.D.

# 1 Commutativity.

**1.0 Definitions**. Let $f{:}A^R \to A$ and $g{:}A^S \to A$ be any two operations on $A$, which again can belong to $\mathcal{O}$ or not. We say that *f and g commute,* when

(14) $$f(g \cdot m) = g(f \cdot \boldsymbol{c}_m), \quad \text{for all } m{:}R \to A^S.$$

[When the carrier is trivial, there only is the empty operation, which must have a nonempty rank and satisfies (14) trivially, since $m \notin (A^S)^R = \emptyset$.] When one of the ranks is empty, say $R = \emptyset$, by (14), (4) and (3) the value of its operation, $f(\emptyset)$, is a "zero" for the other, viz.

(15) $$g(\boldsymbol{k}_{f(\emptyset)}) = f(\emptyset),$$

where now $\boldsymbol{k}_{f(\emptyset)}{:}S \to A$.

In particular, we will be concerned with a *commutative (universal) algebra* that we define by (14), yet now for all $f, g \in \mathcal{O}$. In this algebra, when there is an $R = \emptyset$, (15) holds for every $g$. Then, there is at most one zero. Also, in case of a basis, $X = \emptyset = R$ imply a singleton carrier.

**1.1 Lemmata**.

**(A)** *If f and g commute, then g and f do.*

**(B)** *Given a set $A$, any operation on it commutes with the projections $p_x{:}A^Y \to A$ for all $x \in Y$:*

   (16) $f(p_x \cdot \mathcal{M}) = p_x(f \cdot \boldsymbol{c}_\mathcal{M})$ for all $f{:}A^R \to A$ and $\mathcal{M}{:}R \to A^Y$.

**(C)** *set-ary composition preserves commutativity: given an operation $f{:}A^R \to A$ and the composition $\ell = g \cdot \boldsymbol{c}_G{:}A^Y \to A$ of an indexing $G{:}S \to A^{A^Y}$ of $Y$-ranked operations with an operation $g{:}A^S \to A$, if $f$ commutes with $g$ and with every $G_s$ for $s \in S$, then it does with $\ell$.*

*Proofs.* (A) This symmetry of (14) follows from (6).

(B) [When the carrier is trivial, again we get the self-commutativity of the empty operation as in 1.0.] (B) is trivially true for $Y = \emptyset \not\ni x$. Otherwise and with a nullary $f$, $R = \emptyset = \mathcal{M}$, (16) follows from $f(p_x \cdot \emptyset) = f(\emptyset) \stackrel{(1)}{=} f(\boldsymbol{k}_\emptyset(x)) \stackrel{(0)}{=} (f \cdot \boldsymbol{k}_\emptyset)(x) \stackrel{(9)}{=} p_x(f \cdot \boldsymbol{k}_\emptyset) \stackrel{(4)}{=} p_x(f \cdot \boldsymbol{C}_\emptyset^{(Y)}) = p_x(f \cdot \boldsymbol{c}_\mathcal{M})$.

Otherwise, $[p_x \cdot \mathcal{M}]_r \stackrel{(0)}{=} p_x(\mathcal{M}_r) \stackrel{(9)}{=} \mathcal{M}_r(x) \stackrel{(5)}{=} [\boldsymbol{c}_\mathcal{M}(x)]_r$ for all $r \in R \neq \emptyset$ implies that $p_x \cdot \mathcal{M} = \boldsymbol{c}_\mathcal{M}(x)$, which together with $f(\boldsymbol{c}_\mathcal{M}(x)) \stackrel{(0)}{=} (f \cdot \boldsymbol{c}_\mathcal{M})(x) \stackrel{(9)}{=} p_x(f \cdot \boldsymbol{c}_\mathcal{M})$ gives (16).

(C) [Again, with a trivial carrier self-commutativity occurs because of the remark in 0.4 (A).] When $S = \emptyset$, (15) becomes $f(\boldsymbol{k}_{g(\emptyset)}) = g(\emptyset)$. Then, $f(\ell \cdot \mathcal{M}) \stackrel{(8)}{=} f(\boldsymbol{k}_{g(\emptyset)} \cdot \mathcal{M}) \stackrel{(2)}{=} f(\boldsymbol{k}_{g(\emptyset)}) = g(\emptyset) \stackrel{(1)}{=} \boldsymbol{k}_{g(\emptyset)}(f \cdot \boldsymbol{c}_\mathcal{M}) \stackrel{(8)}{=} \ell(f \cdot \boldsymbol{c}_\mathcal{M})$, for all $\mathcal{M}{:}R \to A^Y$, as required.

When $S \neq \emptyset$, our premises are (14) and

(17) $\qquad f(G_s \cdot \mathcal{M}) = G_s(f \cdot \boldsymbol{c}_{\mathcal{M}}) \quad \text{for all } f{:}A^R \to A,\ \mathcal{M}{:}R \to A^Y \text{ and } s \in S.$

From these premises we have to prove that $f(\ell \cdot \mathcal{M}) = \ell(f \cdot \boldsymbol{c}_{\mathcal{M}})$ for all such $\mathcal{M}$. We consider two cases.

(Case $R = \emptyset$) Notice that, for all $s \in S$, $f(\emptyset) = f(G_s \cdot \emptyset) \stackrel{(17)}{=} G_s(f \cdot \mathbf{C}_\emptyset^{(Y)}) \stackrel{(5)}{=} [\boldsymbol{c}_G(f \cdot \mathbf{C}_\emptyset^{(Y)})]_s$, which by (1) implies $\boldsymbol{c}_G(f \cdot \mathbf{C}_\emptyset^{(Y)}) = \boldsymbol{k}_{f(\emptyset)}$. Then, since $\mathcal{M} = \emptyset$, we get this commutativity as $f(\ell \cdot \emptyset) = f(\emptyset) \stackrel{(15)}{=} g(\boldsymbol{k}_{f(\emptyset)}) = g(\boldsymbol{c}_G(f \cdot \mathbf{C}_\emptyset^{(Y)})) \stackrel{(0)}{=} (g \cdot \boldsymbol{c}_G)(f \cdot \mathbf{C}_\emptyset^{(Y)}) = \ell(f \cdot \mathbf{C}_\emptyset^{(Y)}) \stackrel{(4)}{=} \ell(f \cdot \boldsymbol{c}_\emptyset)$.

(Case $R \neq \emptyset$) For each $\mathcal{M}$, if we define an $m{:}R \to A^S$ by $m_r(s) = G_s(\mathcal{M}_r)$ for all $r \in R$ and $s \in S$, then we get

(18) $\qquad\qquad\qquad\qquad \boldsymbol{c}_G \cdot \mathcal{M} = m,$

from $[(\boldsymbol{c}_G \cdot \mathcal{M})(r)]_s \stackrel{(0)}{=} [\boldsymbol{c}_G(\mathcal{M}_r)]_s \stackrel{(5)}{=} G_s(\mathcal{M}_r)$, while, from $[\boldsymbol{c}_m(s)]r \stackrel{(5)}{=} m_r(s) = G_s(\mathcal{M}_r) \stackrel{(0)}{=} [G_s \cdot \mathcal{M}]_r$, we get $\boldsymbol{c}_m(s) = G_s \cdot \mathcal{M}$, which implies $(f \cdot \boldsymbol{c}_m)(s) \stackrel{(0)}{=} f(\boldsymbol{c}_m(s)) = f(G_s \cdot \mathcal{M}) \stackrel{(17)}{=} G_s(f \cdot \boldsymbol{c}_{\mathcal{M}}) \stackrel{(5)}{=} [\boldsymbol{c}_G(f \cdot \boldsymbol{c}_{\mathcal{M}})]_s$, i.e.

(19) $\qquad\qquad\qquad\qquad f \cdot \boldsymbol{c}_m = \boldsymbol{c}_G(f \cdot \boldsymbol{c}_{\mathcal{M}}).$

Therefore, $f(\ell \cdot \mathcal{M}) = f(g \cdot \boldsymbol{c}_G \cdot \mathcal{M}) \stackrel{(18)}{=} f(g \cdot m) \stackrel{(14)}{=} g(f \cdot \boldsymbol{c}_m) \stackrel{(19)}{=} g(\boldsymbol{c}_G(f \cdot \boldsymbol{c}_{\mathcal{M}})) \stackrel{(0)}{=} (g \cdot \boldsymbol{c}_G)(f \cdot \boldsymbol{c}_{\mathcal{M}}) = \ell(f \cdot \boldsymbol{c}_{\mathcal{M}})$ as required.
<div style="text-align:right">Q.E.D.</div>

**1.2 Theorem.** *Given any set Y, all pairs of Y-ary elementary functions of a commutative (universal) algebra commute.*

*Proof.* When $Y = \emptyset$, we lack projections and any possible $\emptyset$-ary elementary function has to be zero valued, since it comes from possible nullary operations (without them, the statement is trivially true), as we found in **0.4 (B)** and **1.0**. Therefore, for all $f, g \in \mathcal{L}_Y$, we get (14) where $R, S = \emptyset$. Also, when either the algebra is trivial or its carrier is, (14) trivially holds. Then, we consider an $A$ with more than one element and a $Y \neq \emptyset$.

We first prove that any $f \in \mathcal{O}$ commutes with every $\ell \in \mathcal{L}_Y$. After the basis step (16) we only need to consider such an $\ell{:}A^Y \to A$ when it is the composition $\ell = g \cdot \boldsymbol{c}_G$ of an indexing $G{:}S \to \mathcal{L}_Y$ with a $g \in \mathcal{O}$ of rank $S$, where all the $G_s$ commute with $f$. This comes from **1.1 (C)**.

Then, we show that any $j \in \mathcal{L}_Y$ commutes with every $\ell \in \mathcal{L}_Y$. Again, since (16) also holds for $f = j$, we only consider such an $\ell{:}A^Y \to A$ when it is the composition $\ell = g \cdot \boldsymbol{c}_G$ of an indexing $G{:}S \to \mathcal{L}_Y$ with a $g \in \mathcal{O}$ of rank $S$, where $j$ commutes with all the $G_s$. By the first part of this proof and the symmetry in **1.1 (A)** $j$ also has to commute with $g$. Hence, we can again use **1.1 (C)** (with $R = Y$).
<div style="text-align:right">Q.E.D.</div>

**1.3 Corollary.** *The algebra of the conjugate functions of a commutative based algebra is commutative,* viz., with reference to our basis $U{:}X \to A$, given any function $\mathcal{M}{:}X \to A^X$ (now a "square" array),

(20) $\qquad\qquad\qquad \chi_a(\chi_b \cdot \mathcal{M}) = \chi_b(\chi_a \cdot \boldsymbol{c}_{\mathcal{M}}) \quad \text{for all } a, b \in A.$

*Proof.* From **1.2** and **0.8**.
<div style="text-align:right">Q.E.D.</div>

## 2 Endowed dilatation monoids.

2.0 Definitions. We say that an element $d \in A$ of a based algebra as in 0.6 is a *dilatation indicator* when its conjugate function gives us an endomorphism by equalizing its arguments: $\chi_d \cdot \boldsymbol{k}{:}A \to A$ is an endomorphism $\delta \in \mathbb{E}_\alpha = \mathbb{E}_\chi$ for $\alpha$ or for $\chi$ (by 0.8). Then, we also say that $\delta$ is a *dilatation* of our based universal algebra or of $\chi$, while $d$ is *its indicator*.

When $X = \emptyset$, we have a dilatation indicator only when the carrier is singleton, $A = \{d\}$, as in trivial vector spaces with the (zero-valued) dilatation. In fact, for a bigger $A$ and $X = \emptyset$, the functions $\boldsymbol{k}{:}A \to 1$ as in 0.2 (A) and $\chi_d{:}1 \to A$ as in (12) prevent that $\chi_d \cdot \boldsymbol{k} \in \mathbb{E}_\alpha$, since $\mathbb{E}_\alpha = \{\mathbb{I}_A\}$. Yet, we say that $\mathbb{I}_A$ is *the dilatation* for $\chi{:}A \to A^1$, even for a non-singleton $A$, namely even when there are no indicators.

Thus, we are defining (general) dilatations in two different ways depending on the existence of their indicators. Yet, as explained in 2.1 of [16], this split definition comes from the unsplit one of 2.5 of [12] for general universal algebras. Anyway, whatever $X$ is, by 0.8 and 0.4 (B) we easily get that the dilatations exactly are the endomorphisms that also are (rank-less) elementary functions as in 0.4 (C), i.e. their set is $\mathbb{E}_\alpha \cap \mathcal{L}'$.

In case of a vector space, $\mathcal{L}'$ is the set $F$ of all $\delta{:}A \to A$ that multiply the vectors by some scalar of its field, because of the equations that concern the field and the space. Then, the set of vector space dilatations merely is $\mathcal{L}' = F$, since $F \subseteq \mathbb{E}_\alpha$.

Indicators serve to determine the "amount" of a dilatation *with respect to a reference frame* by an element, instead of by a matrix. Yet, while a dilatation by (10) has a single matrix, in general it has a set of indicators, possibly an empty one. $I_\delta \subseteq A$ will denote the set of indicators of dilatation $\delta$ and $\Delta \subseteq \mathbb{E}_\alpha$ will denote the set of all dilatations, while $D \subseteq A$ will denote the set of all dilatation indicators, $D = \bigcup_{\delta \in \Delta} I_\delta$.

Then, we call the *dilatation generator* the function $\gamma{:}D \to \Delta$, such that $\gamma_d = \chi_d \cdot \boldsymbol{k}$ for all $d \in D$ (when $X \neq \emptyset$, $\gamma{:}D \twoheadrightarrow \Delta$). If every $d \in A$ is a dilatation indicator, $D = A$, then we say that the carrier is *dilatation full* and that $\gamma{:}A \twoheadrightarrow \Delta$ is the *total dilatation generator*.

As shown in **1.7** and 2.0 of [16], the most familiar dilatation fullness occurs in vector spaces, where every vector $v \in A$ indicates the dilatation $\gamma_v$ corresponding to the (finite) sum $s$ of its (non-null) components (with respect to the coordinating function for $U$). The value of $\gamma_v$ is the dilatation that multiplies any vector by $s$, see also 3.3.

(One might directly guess this by recalling that $s$ is the linear invariant of the matrix with constant columns $\boldsymbol{k}_v{:}X \to A$, since, when $X$ is a natural number $n+1$, $\gamma_v^{n+1} = s\gamma_v^n$ by the Cayley–Hamilton theorem.)

2.1 Recalled properties.

**(A)** *Under functional composition dilatations form a submonoid of the endomorphism monoid.* (Proof in 2.6 of [16].)

**(B)** *The dilatation monoid is commutative.* (Proof in 2.8 (B) ibid.) Hence,

**(C)** *the dilatation monoid has at most one constant dilatation.* (Proof in 2.9 ibid.)

(Other properties of this monoid concern other generalizations to the universal case of geometrical notions (flocks and scalars), see section **2** of [16] and 1.2–3 of [17].)

2.2 Definitions. In 2.4 and 2.5 we will generalize the construction (**3** in [16]) of the "underlying" field from a vector space, where the field inherits its sum from that of the vector space. We will do it by a minor variant of the following preliminary notions. Then, for each operation $f:A^R \to A$ of our (general) based algebra, we would like to have a corresponding operation $\phi:\Delta^R \to \Delta$ with the same rank, on the set of dilatations.

All such preliminary notions concern *any our based algebra*, yet their terminology comes from vector spaces, which will serve to exemplify them. Then, each general definition will precede a recall of its vector space instance. Such reminders will begin with the words "within our vector space", whereas the general definition with an "in general". In both our based algebra and the vector space $A$ will denote the carrier and $\Delta$ the set of dilatations.

As in [0], we give this (arbitrary) vector space by a sum and all the (rank-less) multiplications times some scalar in its field. Then, according to the convention in 0.3, we denote the sum by an (infix) $+$ and the set of all such rank-less operations by an $F$, $F = \mathcal{O} \cap A^A$, where by the reminder in 2.0 $F$ also is the set of dilatations: $\Delta = F$. Then, within our vector space a $\phi$ will be either another sum, $+:\Delta \times \Delta \to \Delta$, or a rank-less $\phi:\Delta \to \Delta$.

In general, given any $\delta \in \Delta$, let $\boldsymbol{b}_\delta$ denote the function $\boldsymbol{b}_\delta:\Delta \to A^A$ that composes other dilatations with $\delta$: $\boldsymbol{b}_\delta(\beta) = \delta \cdot \beta$ for all $\beta \in \Delta$. Then, $\boldsymbol{b}_\delta:\Delta \to \Delta$ by 2.1 (A). We say that *the dilatation composition distributes over $\phi$*, when, for all $\delta \in \Delta$,

$$(21) \qquad \delta \cdot \phi(\varepsilon) = \phi(\boldsymbol{b}_\delta \cdot \varepsilon), \quad \text{for all } \varepsilon:R \to \Delta,$$

namely when every $\boldsymbol{b}_\delta$ is an endomorphism of $\phi$.

Within our vector space and with a doubleton rank, $R = 2$, $\phi$ corresponds to the latter $+$ through the natural bijection of $\Delta^2 \simeq \Delta \times \Delta$ and this endomorphic property states that the product $\cdot$ of the dilatation monoid distributes over this sum as in $\delta \cdot (\beta + \zeta) = \delta \cdot \beta + \delta \cdot \zeta$ with $\beta, \delta, \zeta \in \Delta$. (Case $R = 1$ is uninteresting, because it restates the commutativity in 2.1 (B).) In general, it will be a *homogeneous* distributivity, as it concerns operations on dilatations only.

Now, let us introduce heterogeneous distributivities. In general, for an $\varepsilon:R \to \Delta$ consider $\boldsymbol{c}_\varepsilon:A \to A^R$ by (7). If

$$(22) \qquad \phi(\varepsilon)(a) = f(\boldsymbol{c}_\varepsilon(a)), \quad \text{for all } \varepsilon:R \to \Delta \text{ and all } a \in A, \text{ and}$$

$$(23) \qquad \delta(f(a)) = f(\delta \cdot a), \quad \text{for all } \delta \in \Delta \text{ and all } a:R \to A,$$

then we say that the application of an element to a dilatation *is (fully) distributive*.

Within our vector space, when $f$ and $\phi$ correspond to (binary) vector sum and to field sum respectively through trivial bijections as above, these two properties provide the scalar times vector product with its two distributive properties, as in $(s + s')v = sv + s'v$ and $s(v + w) = sv + sw$ respectively, where $v, w \in A$ and $s$ and $s'$ are scalars that define two dilatations. (Again, $R = 1$ is uninteresting, because of 2.1 (B).) Hence, we can consider (21), (22) and (23) as generalized distributivities for a possible universal case.

In general now, assume that $A$ is a dilatation full carrier and that the total dilatation generator $\gamma:A \twoheadrightarrow \Delta$ is a homomorphism from every $f:A^R \to A$ onto its $\phi:\Delta^R \to \Delta$. Then, we get an algebra on dilatations that preserves the type of the parent algebra, since it consists of all the $\phi$ that correspond to any $f \in \mathcal{O}$. We call it the *dilatation image of* our based algebra.

Within our vector space, its dilatation image has a dilatation sum and all (rank-less) operations $\phi:\Delta \to \Delta$ such that $\phi = \boldsymbol{b}_\delta$ for some $\delta \in \Delta$, since $\Delta = F$.

In general, in addition to a possible dilatation image, an algebra also has the commutative dilatation monoid in 2.1. Together, they form an algebra on $\Delta$ with a type that extends the one of the parent based algebra to include the two monoid operations. Hence, this new construction changes the species of the starting algebra.

Such an enriched dilatation monoid will also satisfy the equations that the dilatation image inherits from the parent algebra through the homomorphism $\gamma$. When it also satisfies the homogeneous distributivities (21) while the heterogeneous ones in (22) and (23) link it with the starting algebra, we will say that it is *the endowed dilatation monoid of* the based algebra. However, we can simplify such conditions because of the following property.

**2.3 Theorem**. *Whenever a based algebra has a dilatation full carrier, its dilatation generator is a homomorphism onto its dilatation image that gives rise to its endowed dilatation monoid.* (Proved in 4.1 of [16].)

**2.4 Definition**. The endowed dilatation monoid above can have more operations than the ones that are necessary to generate its elementary functions. For instance, for a vector space we do not need the $\boldsymbol{b}_\delta:\Delta\to\Delta$ in 2.2. They are the rank-less elementary functions in 0.4 (C) of the algebra on dilatations that we get by only adding the inherited sum monoid to the dilatation one.

As mentioned in 3.4 (A) this reduct of the endowed dilatation monoid is the "underlying" field of the vector space. This procedure allows this space to be naturally defined as a based universal algebra (without any field). In fact, 3.3 of [16] shows that we get the field (whose scalars are dilatations) by merely calling the application of a vector to a dilatation the "scalar times vector multiplication".

Therefore, we slightly extend the last definition of 2.2 by allowing a dilatation monoid to inherit only some operations. Then, we say that any algebra on dilatations that consists of the dilatation monoid and of some operations of the dilatation image is *an endowed dilatation monoid from* the based algebra. Again, dilatation fullness is enough to achieve such a "monoid".

**2.5 Corollary**. *Whenever a based algebra has a dilatation full carrier, its dilatation generator is a homomorphism onto its dilatation image that gives rise to all endowed dilatation monoids from it.*

*Proof.* The operations of such "monoids" are some of the ones in 2.3. Q.E.D.

**2.6 Theorem**. *All commutative based algebras have endowed dilatation monoids.*

*Proof.* From 1.3 and 2.3 or 2.5 we can prove that (20) implies that, for all $a \in A$, the function $\chi_a \cdot \boldsymbol{k}: A \to A$ in 2.0 is an endomorphism, i.e. that the dilatation generator is total, $\gamma: A \twoheadrightarrow \Delta$. In fact, $\gamma$ will become the required homomorphism. Given any $M: X \to A$, let $\mathcal{M}: X \to A^X$ be the function such that $\mathcal{M}_x = M$ for all $x \in X$. Viz., $\mathcal{M} = \boldsymbol{k}_M$, here with $\boldsymbol{k}: A^X \to (A^X)^X$. Then, $\chi_b \cdot \mathcal{M} \stackrel{(3)}{=} \boldsymbol{k}_{\chi_b(M)}: X \to A$ for all $b \in A$.

We also get $\boldsymbol{c}_\mathcal{M} = \boldsymbol{k} \cdot M$, where $\boldsymbol{k}: A \to A^X$, since $[\boldsymbol{c}_\mathcal{M}(x)]_y \stackrel{(5)}{=} \mathcal{M}_y(x) = M_x \stackrel{(1)}{=} \boldsymbol{k}_{M(x)}(y) \stackrel{(0)}{=} [(\boldsymbol{k} \cdot M)(x)]_y$ for all $x, y \in X$. Therefore, from 1.3 we get $(\chi_a \cdot \boldsymbol{k})(\chi_b(M)) \stackrel{(0)}{=} \chi_a(\boldsymbol{k}_{\chi_b(M)}) = \chi_a(\chi_b \cdot \mathcal{M}) \stackrel{(20)}{=} \chi_b(\chi_a \cdot \boldsymbol{c}_\mathcal{M}) = \chi_b(\chi_a \cdot (\boldsymbol{k} \cdot M)) = \chi_b((\chi_a \cdot \boldsymbol{k}) \cdot M)$ for all $a, b \in A$ and $M: X \to A$, which states that $h = \chi_a \cdot \boldsymbol{k} \in \mathbb{E}_\chi = \mathbb{E}_\alpha$ for all $a \in A$ as in (13). Q.E.D.

## 3 Four examples

**3.0 The complete union semilattices** provide dilatation fullness with a fairly simple introductory example. When they are finite, we can formalize them as the algebras in 0.3 with a nullary operation $\mathbf{0}{:}1 \to A$, such that $\mathbf{0}(\emptyset) = \emptyset$, and the binary union $\cup{:}A \times A \to A$, where the carrier $A$ is the set of all subsets of a finite set, which we denote by $X$. (We need only the finite case here, but the infinite case can also be handled by replacing either $\cup$ with an operation with rank $A$ or both present operations with a single "$\mathbf{P}$-operation" as in [10].) Then, for all $h{:}A \to A$,

(24) $\quad h \in \mathbb{E}_\alpha \quad \text{iff} \quad h(\emptyset) = \emptyset \ \text{ and } \ h(Y \cup Z) = h(Y) \cup h(Z) \ \text{ for all } Y, Z \subseteq X.$

(A) To get a basis as in 0.6, let us choose the frame $U{:}X \to A$ such that $U_x = \{x\}$ for all $x \in X$. We define a function $\eta{:}A^X \to A^A$ by $\eta_M(Y) = \bigcup_{y \in Y} M_y$ for all $M{:}X \to A$ and $Y \subseteq X$. Then, to show (10) we only have to show that *this $\eta$ is the extension function, viz. the inverse of $\boldsymbol{r}_U$*.

*Proof.* We have $[\eta \cdot \boldsymbol{r}_U]_h(Y) \stackrel{(0)}{=} \eta_{h \cdot U}(Y) = \bigcup_{y \in Y}[h \cdot U]_y \stackrel{(0)}{=} \bigcup_{y \in Y} h(U_y) = \bigcup_{y \in Y} h(\{y\}) \stackrel{(24)}{=} h(Y) = [\mathbb{I}_{\mathbb{E}_\alpha}]_h(Y)$ for all $h \in \mathbb{E}_\alpha$ and $Y \in A$, while for all $x \in X$ and $M{:}X \to A$ we have $[\boldsymbol{r}_U \cdot \eta]_M(x) \stackrel{(0)}{=} [\boldsymbol{r}_U(\eta_M)]_x = (\eta_M \cdot U)(x) \stackrel{(0)}{=} \eta_M(U_x) = \eta_M(\{x\}) = \bigcup_{y \in \{x\}} M_y = M_x = [\mathbb{I}_{A^X}]_M(x)$. Q.E.D.

In this analytic representation, we can view every $M = \boldsymbol{r}_U(h){:}X \to A$ as the indexing that provides each *"vertex"* $x \in X$ with $M_x \subseteq X$, an arbitrary subset of its *"sons"*. Namely, the matrices $M$ exactly denote the (simple) directed graphs with vertices in $X$. E.g., in fig.1 we can see the table on the left as the graph on the right.

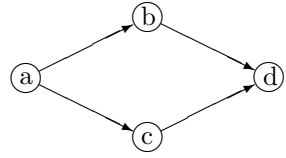

| $x$ | $M_x$ |
|---|---|
| $a$ | $\{b,c\}$ |
| $b$ | $\{d\}$ |
| $c$ | $\{d\}$ |
| $d$ | $\emptyset$ |

Fig. 1

(B) To get the dilatation fullness of $\alpha$ and its endowed dilatation monoid, we use 2.6. In fact, clearly $\alpha$ is commutative. Then, *we get the total dilatation generator as the $\gamma{:}A \twoheadrightarrow \Delta$ such that $\Delta = \{\boldsymbol{k}_\emptyset, \mathbb{I}_A\}$ and $\gamma_Y = \begin{cases} \boldsymbol{k}_\emptyset & \text{when } Y = \emptyset \text{ or} \\ \mathbb{I}_A & \text{when } Y \neq \emptyset \end{cases}$.*

*Proof.* For all $Y, Z \in A$, $\gamma_Y(Z) \stackrel{(0)}{=} \chi_Y(\boldsymbol{k}_Z) \stackrel{(11)}{=} \eta_{\boldsymbol{k}_Z}(Y) \stackrel{(A)}{=} \bigcup_{y \in Y} \boldsymbol{k}_Z(y) \stackrel{(1)}{=} \bigcup_{y \in Y} Z = \begin{cases} \emptyset \stackrel{(1)}{=} \boldsymbol{k}_\emptyset(Z) & \text{when } Y = \emptyset \text{ or} \\ Z = \mathbb{I}_A(Z) & \text{when } Y \neq \emptyset \end{cases}$. Q.E.D.

Therefore, we can use $\gamma$ as in 2.3 to find that, *when $X \neq \emptyset$, the endowed dilatation monoid is a two-elements bounded lattice*.

*Proof.* The unit $\mathbb{I}_A$ of the dilatation monoid and $\boldsymbol{k}_\emptyset$, the $\gamma$-image of the zero of the $\cup$-semilattice, respectively become the unit and zero elements of the lattice. Then, the dilatation composition becomes the meet by (1) and (0), while the $\gamma$-image of $\cup$ becomes the join by the idempotency and the bounds of $\cup$. Q.E.D.

(C) The characteristic function of (sub)sets $j{:}A\!\Vdash\!\twoheadrightarrow\! 2^X$ isomorphically transforms the $\cup$-semilattice in (A) onto another semilattice on $2^X$. Clearly, when we start from the latter semilattice, its analytic representation gives the representation of directed graphs by incidence matrices (up to the usual bijection for $(2^X)^X \simeq 2^{X \times X}$).

We handle such incidence matrices by the usual two-element lattice on $2 = \{0, 1\}$, which is a reduct of the bounded lattice above, up to an obvious isomorphism. This is what occurred in vector spaces: even when we start from an (unusual) vector space as in 0.7 (B), an endowed dilatation monoid as in 2.4 (a field) allows us to build a usual space as in 0.7 (A) that is isomorphic to the original vector space.

**3.1 CPM-PERT networks** ([7] and **3.0** of [13]) show how proper mathematical fittings can enrich the previous directed graphs, while keeping the algebra freedom in 3.0 (A) and the dilatation fullness in 3.0 (B).

We start from a finite set $X$, whose elements are now considered to be *"events"*, and we define an algebra carrier $A$ as the set of all partial functions $a$ from $X$ to the set $N$ of natural numbers, $a{:}Y \to N$ with $Y \subseteq X$. This will replace the sets of "pure" arcs $\langle x, y \rangle$ for $y \in M_x$ in 3.0 with sets of *"activities"*, viz. of arcs $\langle x, y \rangle$ that bear labels that denote *"execution times"* : $t = a(y)$ where $a = M_x$. We call any such time assignment $a \in A$ a *partial schedule* and we denote its domain $Y$ as Dom $a$.

On the carrier $A$ we define three operations: a nullary *zero* **0** as in 3.0, a binary *join* $\sqcup{:}A \times A \to A$ and a (parallel) *successor* $\mathbf{s}{:}A \to A$. The successor is defined for all $a{:}Y \to N$ with $Y \subseteq X$ by the partial schedule $\mathbf{s}(a){:}Y \to N$ such that $[\mathbf{s}(a)]_y = a_y + 1$ for all $y \in Y$. (Then, $\mathbf{s}(\emptyset) = \emptyset$.) The join is defined for all $a{:}Y \to N$ and $b{:}Z \to N$ with $Y, Z \subseteq X$ by the partial schedule $a \sqcup b{:}Y \cup Z \to N$ such that $(a \sqcup b)(x) = \max\{t \in N \mid t = a_x \text{ or } t = b_x\}$ for all $x \in Y \cup Z$. Therefore, now $h{:}A \to A$ is an endomorphism when

(25) $\quad h(\emptyset) = \emptyset, \; h(a \sqcup b) = h(a) \sqcup h(b) \;\text{ and }\; h(\mathbf{s}(a)) = \mathbf{s}(h(a)), \;\text{ for all } a, b \in A.$

To conform to CPM-PERT usage, we use the notations $\max L$ or $\max_{x \in I} \ell_x$ for any finite $L \subseteq N$ or $\ell{:}I \to N$. (According to our set-theoretical model of *natural numbers* [8], they respectively denote $\bigcup L$ or $\bigcup_{x \in I} \ell_x$.) Later we shall use $\cup$ for the binary max.

(A) Set $U_x = \{\langle x, 0 \rangle\}$ for all $x \in X$ to get a frame $U{:}X \to A$ and define $\eta{:}A^X \to A^A$ in agreement with the iteration step of the CPM-PERT algorithm (from **3.0** of [13]), namely by the partial schedules $\eta_M(a){:}\bigcup_{x \in \mathrm{Dom}\, a} D_x \to N$ for all $M{:}X \to A$ and $a \in A$, where $D_x = \mathrm{Dom}\,(M_x)$, such that $[\eta_M(a)]_y = \max\{M_x(y) + a(x) \mid x \in \mathrm{Dom}\, a \text{ and } D_x \ni y\}$ for all $y \in \bigcup_{x \in \mathrm{Dom}\, a} D_x$, Here, $a \in A$ can be thought of as a "preceding partial schedule". Again, *this $\eta$ is the inverse of $\boldsymbol{r}_U$.*

*Proof.* To adapt the proof in 3.0 (A) we first prove three equalities.

Clearly, $\sqcup$ and $\emptyset$ act as the operations of a semilattice whose suprema are computed as the partial schedules $\bigsqcup B{:}\bigcup_{b \in B} \mathrm{Dom}\, b \to N$ such that $[\bigsqcup B]_y = \max\{b_y \mid b \in B \text{ and } \mathrm{Dom}\, b \ni y\}$ for all finite $B \subseteq A$ and $y \in \bigcup_{b \in B} \mathrm{Dom}\, b$. Given $a \in A$ and $x \in \mathrm{Dom}\, a$, let $\lfloor a, x \rfloor$ denote the singleton function $\{\langle x, a(x) \rangle\}$. Then, for all $a \in A$, $a = \bigcup_{x \in \mathrm{Dom}\, a} \lfloor a, x \rfloor = \bigsqcup_{x \in \mathrm{Dom}\, a} \lfloor a, x \rfloor$, since $a' \cup a'' = a' \sqcup a''$ when $\mathrm{Dom}\, a' \cap \mathrm{Dom}\, a'' = \emptyset$. Hence, by the first two equations of (25), $h(a) = \bigsqcup_{x \in \mathrm{Dom}\, a} h(\lfloor a, x \rfloor)$ for all $a \in A$ and $h \in \mathbb{E}_\alpha$, which for all $y \in \mathrm{Dom}\,(h(a))$ implies

(26) $\quad [h(a)]_y = \max\{[h(\lfloor a, x \rfloor)]_y \mid x \in \mathrm{Dom}\, a \text{ and } \mathrm{Dom}\,[h(\lfloor a, x \rfloor)] \ni y\}.$

Notice also that arithmetic induction extends the definition of $\mathbf{s}$ and the last equation of (25) respectively into $[\mathbf{s}^{(n)}(b)]_y = b(y) + n$ and $h(\mathbf{s}^{(n)}(b)) = \mathbf{s}^{(n)}(h(b))$, for all $b \in A$, $y \in \mathrm{Dom}\, b$, $n \in N$ and $h \in \mathbb{E}_\alpha$, where $\mathbf{s}^{(n)}$ denotes the $n$-th composition power of $\mathbf{s}$. Then, for all such $h$, $b$ and $n$,

(27) $\qquad\qquad\qquad \mathbf{s}^{(n)}(U_x) \;=\; \{\langle x, n\rangle\} \qquad$ for all $x \in X$ and

(28) $\qquad\qquad\qquad [h(\mathbf{s}^{(n)}(b))]_y \;=\; [h(b)]_y + n \quad$ for all $y \in \mathrm{Dom}\,(h(b))$.

Now we can check that $\eta$ is the inverse of $\boldsymbol{r}_U$. For all $h \in \mathbb{E}_\alpha$, $a \in A$ and $y \in \mathrm{Dom}\,(h(a))$, $[[\mathbb{I}_{\mathbb{E}_\alpha}]_h(a)]_y = [h(a)]_y \stackrel{(26)}{=} \max\{[h(\lfloor a,x\rfloor)]_y \mid x \in Y \text{ and } D_x \ni y\} \stackrel{(27)}{=} \max\{[h(\mathbf{s}^{(a(x))}(U_x))]_y \mid x \in Y \text{ and } D_x \ni y\} \stackrel{(28)}{=} \max\{[h(U_x)]_y + a(x) \mid x \in Y \text{ and } D_x \ni y\} \stackrel{(0)}{=} \max\{[h \cdot U]_x(y) + a(x) \mid x \in Y \text{ and } D_x \ni y\} = [\eta_{h\cdot U}(a)]_y \stackrel{0.6}{=} [\eta_{\boldsymbol{r}_U(h)}(a)]_y \stackrel{(0)}{=} [[\eta \cdot \boldsymbol{r}_U]_h(a)]_y$, where $Y = \mathrm{Dom}\, a$ and $D_x = \mathrm{Dom}\,[h(\lfloor a,x\rfloor)] = \mathrm{Dom}\,[h \cdot U]_x$ by (27), (28) and (0). Also, for all $M{:}X \to A$, $z \in X$ and $y \in D_z = \bigcup_{x\in\{z\}} \mathrm{Dom}\,(M_x)$, we get $[[\boldsymbol{r}_U \cdot \eta]_M(z)]_y \stackrel{(0)}{=} [\boldsymbol{r}_U(\eta_M)]_z(y) = [(\eta_M \cdot U)(z)]_y \stackrel{(0)}{=} [\eta_M(U_z)]_y = [\eta_M(\{\langle z,0\rangle\})]_y = \max\{M_x(y) + \{\langle z,0\rangle\}(x) \mid x \in \{z\} \text{ and } D_x \ni y\} = M_z(y) + 0 = [[\mathbb{I}_{A^X}]_M(z)]_y$, since $\eta{:}A^X \to \mathbb{E}_\alpha$. $\hfill$ Q.E.D.

In this analytic representation, we can view every $M = \boldsymbol{r}_U(h){:}X \to A$ as the indexing that provides each event $x \in X$ with the indexing $M_x{:}Y \to N$ of the execution times of the activities for the set $Y \subseteq X$ of the events that follow $x$. Namely, the matrices $M$ exactly denote the CPM-PERT projects with the set $X$ of events. E.g., in fig.2 we can see the table on the left as the project on the right.

| $x$ | $M_x$ |
|---|---|
| $a$ | $\{\langle b,1\rangle, \langle c,3\rangle\}$ |
| $b$ | $\{\langle d,7\rangle\}$ |
| $c$ | $\{\langle d,2\rangle\}$ |
| $d$ | $\emptyset$ |

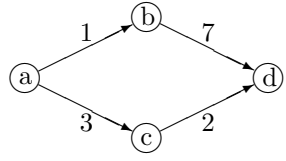

Fig. 2

(B) Again, we use 2.6. We prove that $\alpha$ is *commutative*.

*Proof.* As far as $\mathbf{0}$ and $\sqcup$ are concerned, the proof is straightforward. For $\mathbf{s}$ we get: $\mathbf{s}(\mathbf{s}(a) = \mathbf{s}(\mathbf{s}(a), \mathbf{s}(\mathbf{0}(\emptyset)) = \emptyset = \mathbf{0}(\emptyset)$ as in (15) and $[\mathbf{s}(a \sqcup b)]_x = \max\{t \in N \mid t = a_x \text{ or } t = b_x\} + 1 = \max\{t \in N \mid t = [\mathbf{s}(a)]_x \text{ or } t = [\mathbf{s}(b)]_x\} = [\mathbf{s}(a) \sqcup \mathbf{s}(b)]_x$ for all $a, b \in A$ and $x \in \mathrm{Dom}\, a \cup \mathrm{Dom}\, b$. $\hfill$ Q.E.D.

Then, *we get the total dilatation generator as the* $\gamma{:}A \twoheadrightarrow \Delta$ *such that* $\Delta = \{\boldsymbol{k}_\emptyset\} \cup \bigcup_{n\in N}\{\mathbf{s}^{(n)}\}$ *and* $\gamma_a = \begin{cases} \boldsymbol{k}_\emptyset & \text{when } a = \emptyset \text{ and} \\ \mathbf{s}^{(\mu(a))} & \text{when } a \neq \emptyset \end{cases}$, *where* $\mu(a) = \max_{x \in \mathrm{Dom}\, a} a(x)$, *for all* $a \in A$.

*Proof.* $\gamma_\emptyset(b) \stackrel{(11)}{=} \eta_{\boldsymbol{k}_b}(\mathbf{0}(\emptyset)) \stackrel{(A)}{=} \emptyset \stackrel{(1)}{=} \boldsymbol{k}_\emptyset(b)$ for all $b \in A$ and $[\gamma_a(b)]_y \stackrel{(0)}{=} [\chi_a(\boldsymbol{k}_b)]_y \stackrel{(11)}{=} [\eta_{\boldsymbol{k}_b}(a)]_y \stackrel{(A)}{=} \max\{[\boldsymbol{k}_b(x)]_y + a(x) \mid x \in \mathrm{Dom}\, a \text{ and } D_x \ni y\} \stackrel{(1)}{=} \max\{b(y) + a(x) \mid x \in \mathrm{Dom}\, a \text{ and } \mathrm{Dom}\, b \ni y\} = b(y) + \max_{x \in \mathrm{Dom}\, a} a(x) = [\mathbf{s}^{\mu(a)}(b)]_y$, for all $a, b \in A$ and all $y \in \bigcup_{x \in \mathrm{Dom}\, a} D_x = \bigcup_{x \in \mathrm{Dom}\, a} \mathrm{Dom}\, b = \begin{cases} \emptyset & \text{when } a = \emptyset \text{ and} \\ \mathrm{Dom}\, b & \text{when } a \neq \emptyset \end{cases}$. $\hfill$ Q.E.D.

(C) Every dilatation $\delta \neq \boldsymbol{k}_\emptyset$ defers all execution times of any partial schedule by a constant delay $n$. Then, we could identify it by such an $n$. Yet, if we want to represent all dilatations by natural numbers, then we can also identify it by $n+1$, in order to keep 0 for the empty schedule: $\boldsymbol{k}_\emptyset = \gamma_\emptyset \mapsto 0 = \emptyset$.

Hence, *when $X \neq \emptyset$, the CPM-PERT algebra in (A) has an endowed dilatation monoid that is isomorphic onto the algebra on $N$ that consists of the bounded semilattice of union, of a unary operation $s{:}N \to N$ and of the monoid with the 1-valued nullary constant and a binary operation $\oplus{:}N \times N \to N$, where $s$ and $\oplus$ are respectively defined for all $n, m \in N$ by* $s(n) = \begin{cases} 0 & \text{when } n=0 \\ n+1 & \text{when } n \neq 0 \end{cases}$ *and* $n \oplus m = \begin{cases} n+m-1 & \text{when } n, m \neq 0 \\ 0 & \text{otherwise} \end{cases}$.

*Proof.* If we set $\nu(a) = \begin{cases} 0 & \text{when } a = \emptyset \\ \mu(a)+1 & \text{when } a \neq \emptyset \end{cases}$ for all $a \in A$, then we clearly define a surjection $\nu{:}A \twoheadrightarrow N$, which is the composition of $\gamma$ and of the previous identification map. Let us check that $\nu$ is a homomorphism from $\boldsymbol{0}$ and $\sqcup$ onto the above-claimed bounded semilattice.

$\nu(\boldsymbol{0}(\emptyset)) = \nu(\emptyset) = 0$. Also, since $a \sqcup b = \emptyset$ if and only if $a = b = \emptyset$, we get that $\nu(\emptyset \sqcup \emptyset) = 0 = 0 \cup 0 = \nu(\emptyset) \cup \nu(\emptyset)$, while for all $a, b \in A$ with $Y = \text{Dom } a \cup \text{Dom } b \neq \emptyset$ we still get $\nu(a \sqcup b) = \max_{x \in Y}(\max\{t \in N \mid t = a_x \text{ or } t = b_x\}) + 1 = ((\max_{x \in \text{Dom } a} a_x) \cup (\max_{x \in \text{Dom } b} b_x)) + 1 = (\mu(a) \cup \mu(b)) + 1 = (\mu(a)+1) \cup (\mu(b)+1)$, where $\mu(a)+1 = \nu(a) > 0$ or $\mu(b)+1 = \nu(b) > 0$. Hence, $\nu(a \sqcup b) = \nu(a) \cup \nu(b)$.

Also, $\nu$ is a homomorphism from $\boldsymbol{s}$ onto $s$, $\nu(\boldsymbol{s}(a)) = s(\nu(a))$ for all $a \in A$, because $\nu(\boldsymbol{s}(a)) = \begin{cases} 0 = s(0) = s(\nu(a)) & \text{when } a = \emptyset \\ \mu(\boldsymbol{s}(a))+1 = \mu(a)+2 = \nu(a)+1 = s(\nu(a)) & \text{when } a \neq \emptyset \end{cases}$.

Finally, to check that 1 and $\oplus$ identify the dilatation monoid, consider the identification map and the identities: $\boldsymbol{s}^{(n)} \cdot \boldsymbol{s}^{(m)} = \boldsymbol{s}^{(n+m)}$ and $(\boldsymbol{s}^{(n)} \cdot \boldsymbol{k}_\emptyset)(a) \stackrel{2.1(B)}{=} (\boldsymbol{k}_\emptyset \cdot \boldsymbol{s}^{(n)})(a) \stackrel{(0)}{=} \boldsymbol{k}_\emptyset(\boldsymbol{s}^{(n)}(a)) \stackrel{(1)}{=} \emptyset \stackrel{(1)}{=} \boldsymbol{k}_\emptyset(\boldsymbol{k}_\emptyset(a)) \stackrel{(0)}{=} (\boldsymbol{k}_\emptyset \cdot \boldsymbol{k}_\emptyset)(a)$ for all $n, m \in N$ and $a \in A$, which give the required monoid identities through $\nu$.  Q.E.D.

This representation of $\Delta$ by $N$ gives a natural bounded semilattice, yet the other operations (the nullary one, $s$ and $\oplus$) are a bit clumsy. One might well replace $N$ with its successor set $\mathcal{S}(N) = N \cup \{N\}$, to identify every $\boldsymbol{s}^{(n)}$ by the "delay" $n \in N$ and $\boldsymbol{k}_\emptyset$ by $N$, the "undefinable delay". Then, the latter operations become natural (within cardinal arithmetics), yet the lattice ones now become clumsy.

(D) Again, *an isomorphism $j$ transforms our algebra on $A$ onto another on $\Delta^X$*.

*Proof.* Take $j_a(x) = \begin{cases} \boldsymbol{k}_\emptyset & \text{if } x \notin \text{Dom } a \\ \boldsymbol{s}^{(a(x))} & \text{otherwise} \end{cases}$ for all $a \in A$ and $x \in X$. Then, $j{:}A \twoheadrightarrow \Delta^X$, since $j$ has the inverse $\ell{:}\Delta^X \twoheadrightarrow A$ such that $\text{Dom } \ell_\delta = \{x \in X \mid \delta(x) \neq \boldsymbol{k}_\emptyset\}$ and $\boldsymbol{s}^{(\ell_\delta(x))} = \delta(x)$ for all $\delta{:}X \to \Delta$ and $x \in \text{Dom } \ell_\delta$.  Q.E.D.

**3.2 The ring of integers** is a simple example of a ring that will come as an endowed dilatation monoid, but not from a module, i.e. not as an underlying ring. Hence, it strengthens the difference between our geometric generalizations and the algebraic ones, by providing the counterexample in [16] about the semi-ring of natural numbers with a converse.

(A) Let $\alpha$ be the group of the integers under addition and consider its extension function. A singleton $U{:}X \to A$ with value 1 is a frame. Then, by the natural map $j'$ for $A^X \simeq A$ we can replace the endomorphism representability (10) with $r{:}\mathbb{E}_\alpha \twoheadrightarrow A$, where $r(h) = h(1)$ for all $h \in \mathbb{E}_\alpha$, and we consider the $\epsilon{:}A \to A^A$ such that $\epsilon_a(b) = ab$ for all $a, b \in A$. Then, to show that $U$ is a reference frame, we only have to show that *this $\epsilon$ replaces the extension* *function, viz. it is the inverse of $r$.*

*Proof.* In fact, $[\epsilon \cdot r]_h(b) \stackrel{(0)}{=} h(1)b = bh(1) = \bar{b}\sum_{i=1}^{i=|b|} h(1) = h(\bar{b}\sum_{i=1}^{i=|b|} 1) = h(b) = [\mathbb{I}_{\mathbb{E}_\alpha}]_h(b)$ for all $h \in \mathbb{E}_\alpha$ and $b \in A$, where $\bar{b}$ and $|b|$ respectively denote the sign of $b$ and the natural number corresponding to its absolute value. Conversely, $[r \cdot \epsilon]_a \stackrel{(0)}{=} \epsilon_a(1) = a1 = a = [\mathbb{I}_A]_a$ for all $a \in A$.  Q.E.D.

(B) Therefore, by 2.6 every integer $a \in A$ gives us a dilatation $\gamma_a \in \Delta$ that corresponds to $\chi_a$ up to $j'$, since $X$ is singleton. Hence, $\gamma = \epsilon$, since $\gamma_a(b) \stackrel{(11)}{=} \epsilon_b(a) = ba = ab = \epsilon_a(b)$ for all $a, b \in A$. This implies that $\Delta = \mathbb{E}_\alpha$ is the set of multipliers by an integer.

Then, the further isomorphism $\epsilon^{-1} = r:\Delta \mapsto\!\!\!\!\!\twoheadrightarrow A$, which is the inverse of $\gamma$, gives us back the integers (with their sum group). Clearly, under $\epsilon^{-1}$ the dilatation monoid becomes the multiplication monoid. To sum up: *the endowed dilatation monoid of the group of integers under addition is their ring, up to an isomorphism.*

(C) Our $\mathbb{E}_\alpha$ is the set of endomorphisms also of the monoid of integers under addition and the representability proof in (A) holds even for this monoid. Yet, $U$ is not a monoid generator. Then, endomorphism representability does not imply carrier generation, contrary to the claim of [11] (its flawed proof in 6.7 ibid. identified two different extension functions).

**3.3 Gaussian integers.** The ring of integers could not give us a meaningful isomorphism $j:A \mapsto\!\!\!\!\!\twoheadrightarrow \Delta^X$, since $X$ was singleton. However, we can, if we start from the group of Gaussian integers under addition. For every $a \in A$ let $a = a' + a''\imath$, where $a'$ and $a''$ are (real) integers. Here, we omit the (usual) quantifications.

Take $X = 2$ with $U_0 = 1$, $U_1 = \imath$ and $\eta_M(a) = M'_0 a' + M'_1 a'' + (M''_0 a' + M''_1 a'')\imath$. Then, (10) comes from $[(\eta \cdot \boldsymbol{r}_U)(h)]_a = \eta_{h \cdot U}(a) = h(U_0)'a' + h(U_1)'a'' + (h(U_0)''a' + h(U_1)''a'')\imath = h(1)'a' + h(\imath)'a'' + (h(1)''a' + h(\imath)''a'')\imath = h(a')' + h(a''\imath)' + (h(a')'' + h(a''\imath)'')\imath = h(a)' + h(a)''\imath = h(a)$, while $(\boldsymbol{r}_U \cdot \eta)_M(0) = M'_0 + M'_1 0 + (M''_0 + M''_1 0)\imath = M'_0 + M''_0 \imath = M_0$ and $(\boldsymbol{r}_U \cdot \eta)_M(1) = M'_0 0 + M'_1 + (M''_0 0 + M''_1)\imath = M_1$.

Also, $\gamma_a(b) = \eta_{\boldsymbol{k}(b)}(a) = b'a' + b'a'' + (b''a' + b''a'')\imath = (a' + a''\imath)(b' + b''\imath) = (a' + a''\imath)b$. Namely, the dilatations are integer multipliers as in 3.2 and we easily get that the previous ring is isomorphic to their endowed monoid. Yet, this time we get a meaningful $j:A\mapsto\!\!\!\!\!\twoheadrightarrow \Delta^2$ by $j(m_0(1) + m_1(\imath)) = m = (m_0, m_1)$  for all $m \in \Delta^2$.

**3.4 The dilatation generator** tells us how many dilatations we have and how their indicators are spread in a carrier. Even though the dilatation idea comes from the Geometry of vector spaces, which are commutative, this generator provides all based universal algebras with useful information. Let us review three cases.

(A) When the carrier is dilatation full, as for our commutative algebras, the ensuing homomorphism in 2.5 always gives an endowed dilatation monoid. Then, it gives us a phylogenic relationship between two algebras with different types.

In [17] this allowed vector spaces to have a natural characterization as universal algebras, which embodies their well-known projective features: *they merely become "dilatation completed" Abelian groups with a dilatable reference frame.* Neither fields nor all equations are necessary, since they stem from this as endowed dilatation monoids, which by 2.5 do not need the equations mentioned in 2.2 (including the ones due to the homomorphism $\gamma$).

The endowed dilatation monoids from dilatation full algebras clearly have (true) polynomial functions, like the ones of fields. Then, among several open problems, one concerns the Cayley–Hamilton theorem. Generalizations so far of this theorem work in all universal algebras [11], but their characteristic equations relate an algebra endomorphism with its "eigen/inner-spaces" without polynomials. (Such equations are made by certain "forms of higher degree" in **4** ibid.) Thus, one might check whether in dilatation full algebras we can also use such polynomials as we conveniently do for characteristic equations in vector spaces.

(In the past century some authors, e.g. [5], used the adjective "polynomial" instead of "elementary" to denote what in a vector space merely is a *linear* term or its function and not what concerns the "underlying" field. This is one of the instances where structures of vector spaces were implicitly denied to belong to universal algebras. Some of such misunderstandings (see 2.1 and (C)) are due to missing developments in Universal Algebra. Yet, they likely began with a misinterpretation of Segre's semi-linear transformations [21]: defining vector spaces as modules [0] allowed (couple of) automorphisms to define invariance, instead of acknowledging that automorphisms fail it even in mere universal algebras, see 0.6.)

(B) When the algebra is "dilatation poor", viz. there only are trivial dilatations, there are also only trivial isomorphisms between dilatation monoids. Then, such algebras have easy transformations, viz. one does not need to choose the scalar isomorphism for their "Segre descriptions" (3.3 in [19]), which generalize Segre's semi-linear transformations.

Two well-known poorness examples concern Computer Science: Boolean algebras and word catenation monoids, which respectively have one and two dilatations, see 2.7 of [16]. Two-ary integers [15] (which model word addressing and Web transformations) form a lesser-known example, in spite of its isomorphism to an early well-known algebra [6].

(C) The third case is intermediate between the former two and overlaps them. It concerns the dilatation indicators of the identity, called *flock combiners* in **2.0** of [18], which, together with their extension in [20], serve to define flocks in every based universal algebra. For instance, in the Boolean case (2.4 (A) ibid. and 3.5 ibid.) the identity has enough flock combiners to give an affine structure similar to the one of affine hulls of vector spaces.

Anyway, such flocks always enjoy the main properties of the well-known flocks of vector spaces (ibid.). Therefore, contrary to the Geometry developments of the past century, one should not restrict the study of affine structures to seemingly immediate generalizations of vector spaces like modules. Universal algebras can work and, contrary to modules from skew fields, exist even in the finite. Furthermore, while they naturally contain vector spaces, they are elementary algebraic structures, not composed ones like modules.

Acknowledgments. The organization of AAA75 allowed the author (who was at the University of Parma then) to present an outline of this work at the Darmstadt conference of November 2–4, 2007.

Gabriele Ricci
*gabrieleericci@bluewin.ch,   gabriele.ricci@unipr.it*
Home phones: (+41/0) 79 325 6125; (+33/0) 4 9316 9423